\documentclass[letterpaper, 10pt, onecolumn]{ieeeconf}

\IEEEoverridecommandlockouts        % This command is only
\usepackage{hyperref}

\usepackage{fancyhdr}
\usepackage{graphicx,color}
\usepackage{amsmath}
\usepackage{amssymb}
\usepackage{mathrsfs}
\usepackage{subfigure}
\usepackage{url}
\usepackage{color}
\usepackage{algorithmic}
\usepackage{dsfont}
\usepackage{bbm}

\usepackage{array}
\usepackage[table]{xcolor}
\usepackage{yfonts}

\newtheorem{example}{Example}

\newcommand{\map}[3]{#1: #2 \rightarrow #3}
\newcommand{\setdef}[2]{\{#1 \; : \; #2\}}

\newcommand{\supscr}[2]{{#1}^{\textup{#2}}}
\newcommand{\until}[1]{\{1,\dots,#1\}}

\newcommand{\real}{\mathbb{R}}

\newcommand{\transpose}{\mathsf{T}} %or \top or \intercal

\newcommand{\mc}{\mathcal}

\DeclareSymbolFont{bbold}{U}{bbold}{m}{n}
\DeclareSymbolFontAlphabet{\mathbbold}{bbold}

\newcommand\oprocendsymbol{\hbox{$\square$}}
\newcommand\oprocend{\relax\ifmmode\else\unskip\hfill\fi\oprocendsymbol}

\usepackage{tikz}
\usetikzlibrary{matrix}

\renewcommand{\theenumi}{(\roman{enumi})}

\title{\LARGE \bf On Structural Controllability of Symmetric (Brain)
  Networks}

\author{Tommaso Menara, Shi Gu, Danielle S. Bassett, and Fabio
  Pasqualetti%
  \thanks{Tommaso Menara and Fabio Pasqualetti are with the Mechanical
    Engineering Department, University of California at Riverside,
    \{\href{mailto:tomenara@engr.ucr.edu}{\texttt{tomenara}},
    \href{mailto:fabiopas@engr.ucr.edu}{\texttt{fabiopas\}@engr.ucr.edu}}.
    Shi Gu is with the Department of Computer Science and Engineering, University of Electronic Science and Technology, Chengdu, China. Danielle S. Bassett is with the Department of Bioengineering, the Department of Electrical and Systems Engineering, the Department of Physics and Astronomy, and the Department of Neurology, University of Pennsylvania,
    \{\href{mailto:dsb@seas.upenn.edu}{\texttt{dsb@seas.upenn.edu}\}}.}}

\begin{document}
\maketitle
\thispagestyle{fancy}
\pagestyle{fancy}

\begin{abstract}
  \noindent
  The question of controllability of natural and man-made network
  systems has recently received considerable attention. In the context
  of the human brain, the study of controllability may not only shed
  light into the organization and function of different neural
  circuits, but also inform the design and implementation of minimally
  invasive yet effective intervention protocols to treat neurological
  disorders. While the characterization of brain controllability is
  still in its infancy, some results have recently appeared and given
  rise to scientific debate. Among these,
  \cite{SG-FP-MC-QKT-BYA-AEK-JDM-JMV-MBM-STG-DSB:15} has numerically
  shown that a class of brain networks constructed from DSI/DTI
  imaging data are controllable from one brain region. That is, a
  single brain region is theoretically capable of moving the whole
  brain network towards any desired target state. In this note we
  provide evidence supporting controllability of brain networks from a
  single region as discussed in
  \cite{SG-FP-MC-QKT-BYA-AEK-JDM-JMV-MBM-STG-DSB:15}, thus
  contradicting the main conclusion and methods developed
  in~\cite{CT-RPR-MC-SZ-MZ-SS:17}.
\end{abstract}

We consider brain networks modeled by a weighted graph
$\mc G = (\mc V, \mc E)$, where $\mc V = \until{n}$ and
$\mc E \subseteq \mc V \times \mc V$ are the vertex and edge sets,
respectively. Let $A = [a_{ij}]$ be the weighted adjacency matrix of
$\mc G$, where $a_{ij} = 0$ if $(i,j) \not\in \mc E$ and
$a_{ij} \in \real_{\ge 0}$ if $(i,j) \in \mc E$. We assume that $A$ is
symmetric and that the graph $\mc G$ has no self loops, forcing the
diagonal entries of $A$ to zero. These assumptions are dictated by the
use of DSI/DTI scans to reconstruct brain networks
\cite{SG-FP-MC-QKT-BYA-AEK-JDM-JMV-MBM-STG-DSB:15}. Let
$\map{x}{\mathbb{N}}{\real^n}$ be the vector containing the state of
the brain regions over time. The network dynamics with control region
$i \in \mc V$ read as
\begin{align}\label{eq: network}
  x(t+1) = A x(t) + \supscr{b}{i} u(t),
\end{align}
where $\map{u}{\mathbb{N}}{\real}$ is the control input injected into
the $i$-th brain region, and the input vector $\supscr{b}{i}$
satisfies $b_j^\text{i} = 0$ if $j \neq i$ and $b_i^\text{i} = 1$. The
network \eqref{eq: network} is controllable if and only if the
controllability matrix $\mc C (A, \supscr{b}{i})$ is invertible
\cite{tk:80}, where
\begin{align}\label{eq: contr_matr}
    \mc C (A, \supscr{b}{i}) &=
    \begin{bmatrix}
      \supscr{b}{i} & A \supscr{b}{i} & \cdots & A^{n-1} \supscr{b}{i}
    \end{bmatrix}
                         .
\end{align}

Assessing controllability of network systems is numerically difficult
because the controllability matrix typically becomes ill-conditioned
as the network cardinality increases; e.g., see
\cite{FP-SZ-FB:13q,js-aem:13}. Because different controllability tests
suffer similar numerical difficulties, a convenient tool to study
controllability of networks is to resort to the theory of structural
systems. To formalize this discussion, notice that the determinant
$\det (\mc C (A, \supscr{b}{i})) = \phi (a_{ij})$ is a polynomial
function of the nonzero entries of the adjacency matrix. The network
\eqref{eq: network} is uncontrollable when the weights are chosen so
that $\mc C (A, \supscr{b}{i})$ is not invertible or, equivalently,
when $\phi (a_{ij}) = 0$. Let $\mc S$ contain the choices of weights
that render the network \eqref{eq: network} uncontrollable, that is,
\begin{align}\label{eq: determinant}
  \mc S = \setdef{a_{ij}}{(i,j)\in \mc E \text{ and }\phi (a_{ij}) = 0}.
\end{align}
Notice that each element of $\mc S$ can be represented as a point in
$\real^d$, where $d = |\mc E|$ is the number of nonzero entries of
$A$. Formally, the set $\mc S$ defines an algebraic variety of
$\real^d$ \cite{wmw:85}. This implies that controllability of
\eqref{eq: network} is a \emph{generic property} because it fails to
hold on an algebraic variety of the parameter space
\cite{mumford1999red}. Thus, when assessing controllability of
\eqref{eq: network} as a function of the network weights, only two
mutually exclusive cases are possible: \smallskip
\begin{enumerate}
\item either there is \emph{no} choice of weights $a_{ij}$, with
  $(i,j)\in \mc E$, rendering the network \eqref{eq: network}
  controllable, or
  \smallskip
\item the network \eqref{eq: network} is controllable for all choices
  of weights $a_{ij}$ except, possibly, those lying in a proper
  algebraic variety of the parameter space $\real^d$ (see Example
  \ref{ex: grid} below).
\end{enumerate}
\smallskip Loosely speaking the above discussion implies that, if one
can find a choice of weights $a_{ij}$ such that \eqref{eq: network} is
controllable, then \emph{almost all} choices of weights $a_{ij}$ yield
a controllable network. In this case, the network is said to be
structurally controllable \cite{wmw:85,lin1974structural,kjr:88}. In
what follows we show that brain networks are structurally controllable
from one single region, thus providing theoretically-validated and
numerically-reliable support to the result in
\cite{SG-FP-MC-QKT-BYA-AEK-JDM-JMV-MBM-STG-DSB:15}. This further shows
that the result in \cite{CT-RPR-MC-SZ-MZ-SS:17} is likely incorrect
and misleading. In fact, even in an unfortunate choice of weights that
prevents controllability, that is, a choice of weights that lies in a
proper algebraic variety, a random and arbitrarily small deviation of
network weights due to perturbation or uncertainty in estimating
neural connections would guarantee controllability.

Classic results on structural controllability cannot be directly
applied to symmetric (brain) networks. In fact, these results assume
that the network weights can be selected arbitrarily and independently
from one another, a condition that cannot be satisfied when the
weights need to be symmetric. To overcome this limitation we proceed
as follows: first we show that network controllability remains a
generic property when the weights are symmetric; then, we find a
choice of symmetric weights that guarantees controllability. This
ensures that brain networks are structurally controllable from one
single node, even with symmetric weights, and that almost all choices
of symmetric edge weights yield controllable networks.

\smallskip
\noindent
\textit{(Generic controllability of networks with symmetric weights)}
Let $d = |\mc E|$ and notice that a network with symmetric weights is
uniquely specified by $d/2$ parameters, which can be represented as a
point in the Euclidean space $\real^{d/2}$. With the symmetry
constraint, the determinant of the controllability matrix is a
polynomial function of $d/2$ parameters and can be obtained, for
instance, from the determinant $\phi (a_{ij})$ in \eqref{eq:
  determinant} by substituting $a_{ij}$ with $a_{ji}$ whenever
$i > j$. Thus, even for symmetric networks, the determinant of the
controllability matrix is a polynomial function, and the weights that
render the network uncontrollable lie on an algebraic variety of the
parameter space $\real^{d/2}$. We conclude that controllability of
networks with symmetric weights remains a generic property: either
there exists no choice of symmetric weights $a_{ij}$ that makes the
network controllable, or the network is controllable for almost all
choices of symmetric weights.

\smallskip
\noindent
\textit{(A controllable realization of a brain network)} Because
controllability of symmetric networks is a generic property, it is
sufficient to construct a controllable symmetric network to show that
almost all choices of weights yield a controllable network. To do
this, we construct a Hamiltonian path starting from the control
node,\footnote{A path in a graph is Hamiltonian if it visits all the
  vertices exactly once.} select the weights of the edges in the path
equal to one, and set all other weights equal to zero. Notice that the
determinant of the controllability matrix associated with the
constructed network has unit magnitude, proving that the network is
structurally controllable with symmetric weights.

\begin{figure}[t]
  \centering 
  \subfigure[]{
    \includegraphics[width=.32\columnwidth]{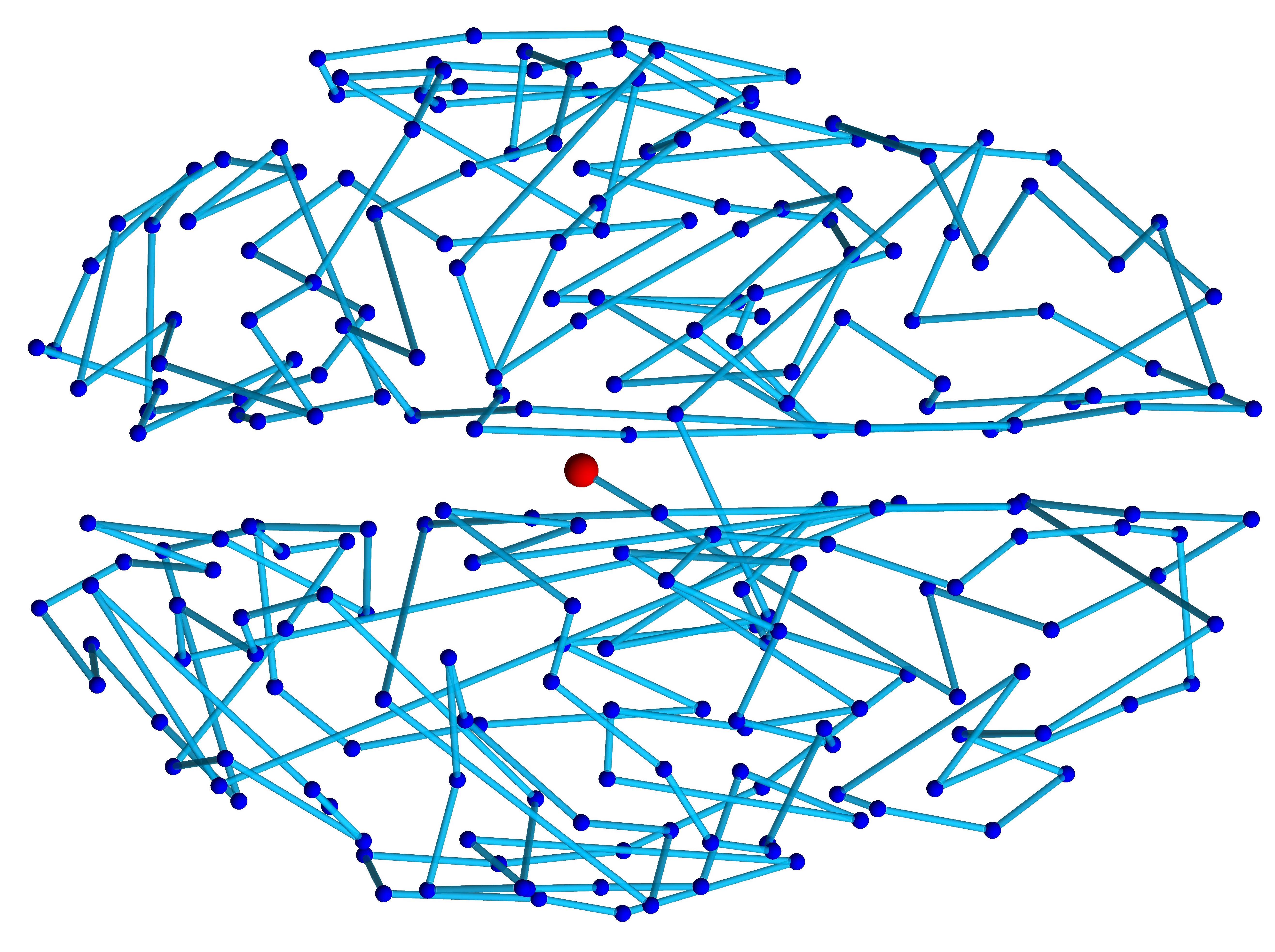}
  }\hfill
    \subfigure[]{
    \includegraphics[width=.25\columnwidth]{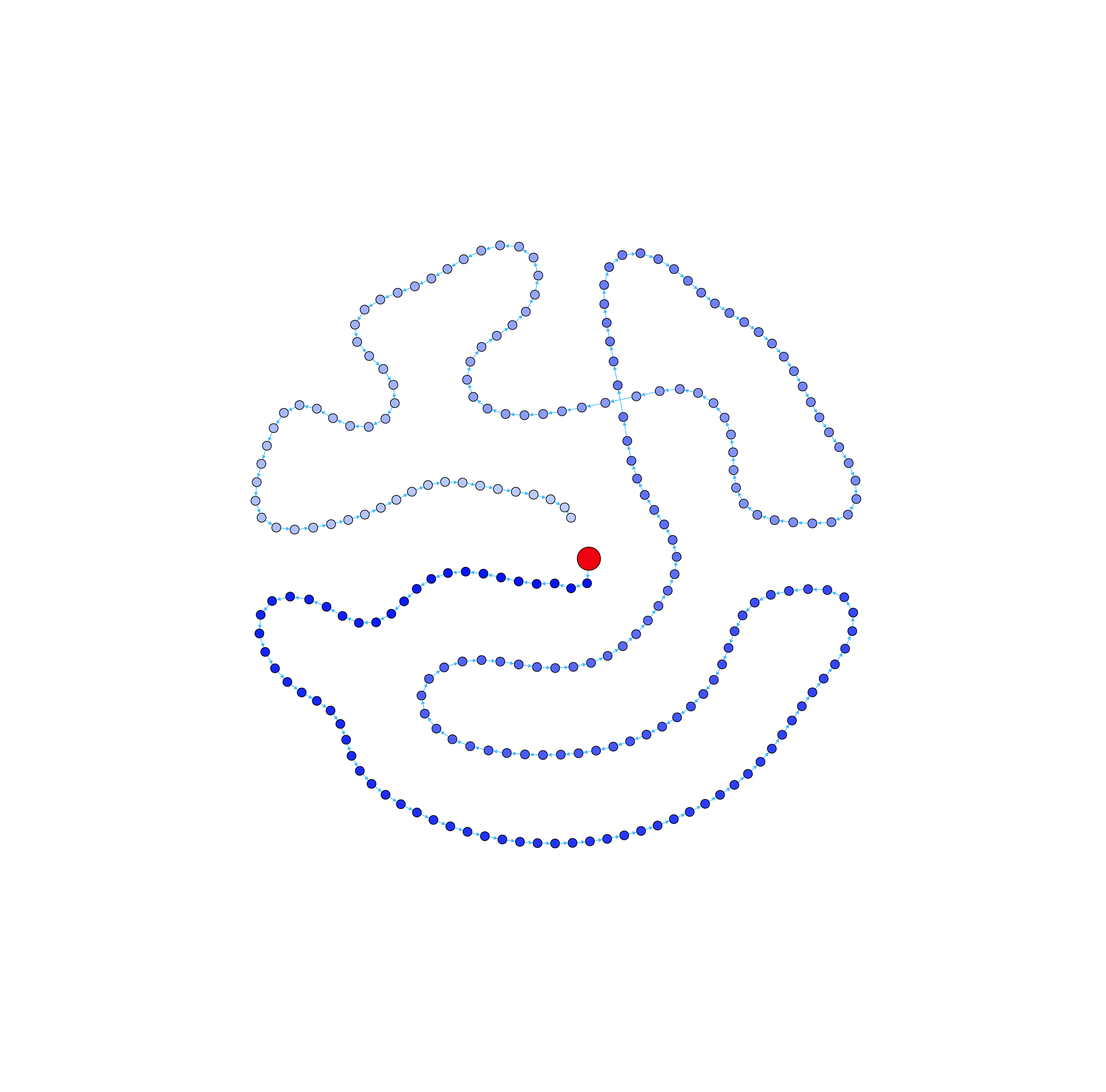}
  } \hfill
  \subfigure[]{
    \includegraphics[width=.25\columnwidth]{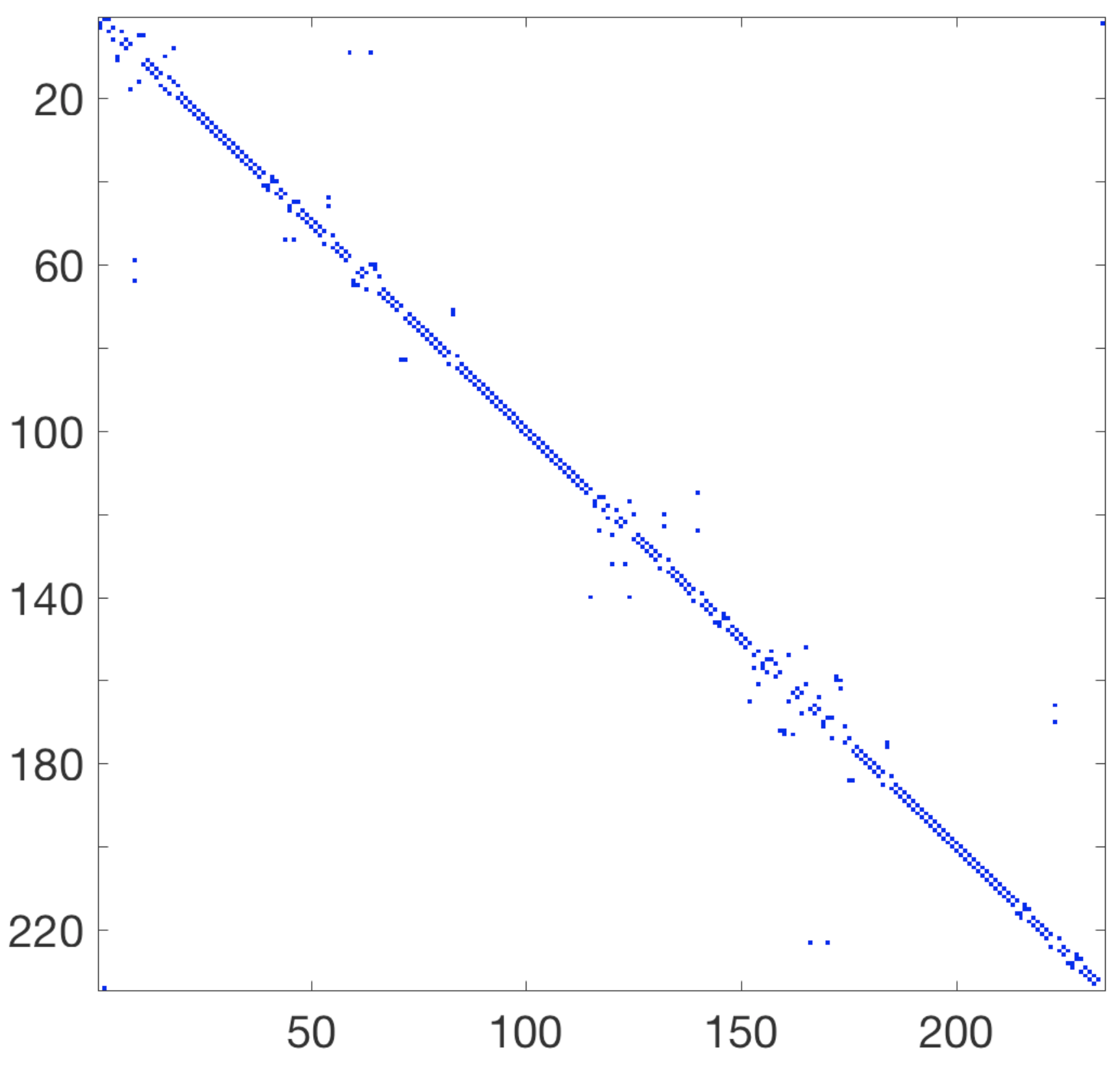}
  } 
  \caption{(a) Hamiltonian path in a structural brain network with 234
    nodes. The control node is highlighted in red. (b) A two
    dimensional representation of the Hamiltonian path. (c) Adjacency
    matrix of the controllable network corresponding to the
    Hamiltonian path with unit weights.}
  \label{fig: hamPaths}
\end{figure}

\smallskip We remark that not all networks admit a Hamiltonian path,
and that the existence of such a path is only a sufficient condition
for structural controllability with symmetric weights. Yet, as we show
in Fig. \ref{fig: hamPaths} for one brain network and one control
region, all the networks in our dataset admit a Hamiltonian path from
every region, showing that they are structurally controllable from
every region even with symmetric weights.  We conclude this discussion
with an academic example.

\begin{example}{\textbf{\emph{(Structural controllability of a
        network with symmetric weights)}}}\label{ex: grid}
  Consider a network with adjacency matrix
  \begin{align}
    A = 
    \begin{bmatrix}
      0 & a_{12} & a_{13}\\
      a_{12} & 0 & a_{23}\\
      a_{13} & a_{23} & 0
    \end{bmatrix},
  \end{align}
  control node $\{1\}$, and input vector $\supscr{b}{1} =
    \begin{bmatrix}
      1 & 0 & 0\\
    \end{bmatrix}^\transpose$. The network is represented in
    Fig. \ref{fig: grid_plain}. Following our analysis, the network is
    structurally controllable even with symmetric weights because it
    admits a Hamiltonian path (see Fig. \ref{fig: grid_path}). That
    is, the network is controllable for almost all choices of
    symmetric weights. To see this, compute the controllability matrix
    \begin{align}
      \mc C (A, \supscr{b}{1}) = 
      \begin{bmatrix}
        1 & 0 & a_{12}^2 + a_{13}^2\\
        0 & a_{12} & a_{13}a_{23}\\
        0 & a_{13} & a_{12}a_{23}
      \end{bmatrix},
    \end{align}
    and the determinant
    $\det (\mc C (A, \supscr{b}{1})) = a_{23} a_{12}^2 - a_{23}
    a_{13}^2$. Thus, the network is controllable (i.e.,
    $\det (\mc C (A, \supscr{b}{1})) \neq 0$) for all symmetric
    choices of weights $a_{12}$, $a_{13}$, and $a_{23}$, except those
    lying on the proper algebraic variety shown in Fig. \ref{fig:
      hypersurface}. \oprocend
\end{example}

\begin{figure}[htb]
  \centering \subfigure[]{
    \includegraphics[width=.22\columnwidth]{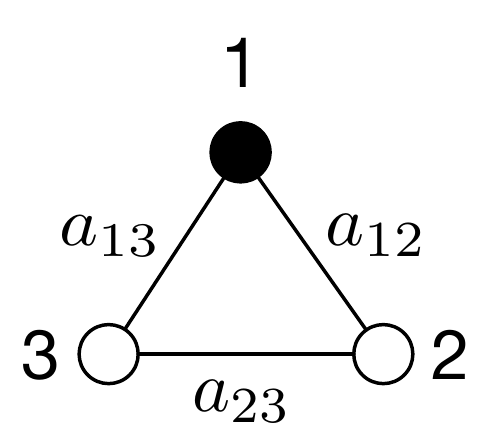} \label{fig:
      grid_plain}
  }\quad\;\;\;\;\;\;\;\;\;\; \;\;\;\;\;\;\;\;\;\;
  \subfigure[]{
    \includegraphics[width=.22\columnwidth]{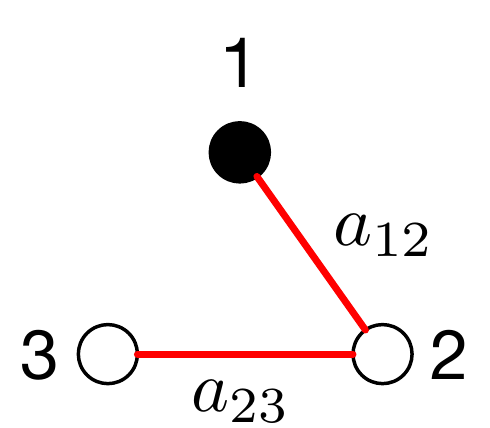}\label{fig:
      grid_path}
  }\hfill
    \subfigure[]{
      \includegraphics[width=.5\columnwidth]{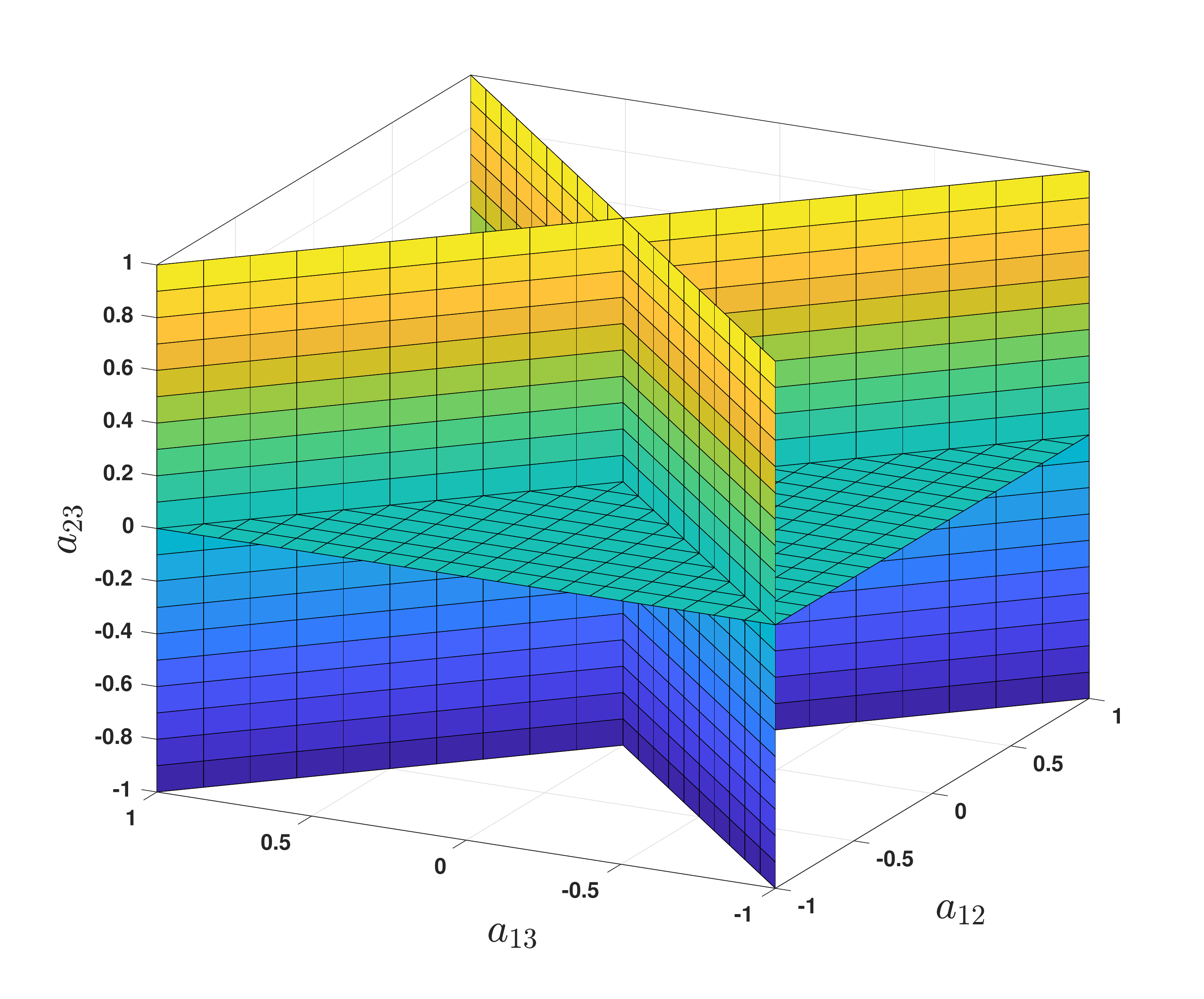}\label{fig:
        hypersurface}
  }
  \caption[Optional caption for list of figures]{(a) Network with
    symmetric weights considered in Example \ref{ex: grid}. (b)
    Network induced by a Hamiltonian path starting from the control
    node. (c) Algebraic variety containing the weights for which the
    network is not controllable. The network is controllable for all
    weights outside of this hypersurface.}
\end{figure}

\bibliographystyle{unsrt}
\bibliography{./BIB}

\end{document}